\documentclass[a4paper,twoside,10pt]{article}

\usepackage{amssymb}
\usepackage{amsmath}
\usepackage{amscd}
\usepackage[matrix]{xypic}
\usepackage{cite}

\newtheorem{theorem}{Theorem}

\newtheorem{definition}[theorem]{Definition}

\newtheorem{proposition}[theorem]{Proposition}

\newtheorem{defi}[theorem]{D\'efinition}

\newcommand\K{\mathbb{K}}

\begin{document}

\noindent {\bf \LARGE{Cog\`ebres Lie-admissibles}}

\bigskip

\bigskip

\noindent Michel Goze\thanks{Email address: M.Goze@uha.fr}, Elisabeth Remm\thanks{Email address: E.Remm@uha.fr}. \ \ 
Laboratoire de Math\'ematiques et Applications, \
Universit\'{e} de Haute Alsace,\
4, rue des Fr\`eres Lumi\`ere, F-68093 Mulhouse Cedex, France.

\bigskip

\noindent {\small  {\bf Mots cl\'es}: Alg\`ebres non associatives, alg\`ebres Lie-admissibles, 
cog\`ebres non associatives.}

\bigskip

\noindent {\bf R\'esum\'e.} Nous d\'efinissons les cog\`ebres associ\'ees 
\`a des classes d'alg\`ebres non associatives dont l'associateur v\'erifie 
des conditions d'invariance donn\'ees par l'action du groupe sym\'etrique d'ordre 3. 
Parmi ces alg\`ebres
nous retrouvons les alg\`ebres de Vinberg et pr\'e-Lie, les alg\`ebres Lie-admissible
et les alg\`ebres \`a puissance 3-associative. Nous \'etudions les propri\'et\'es
liant ces cog\`ebres et alg\`ebres non associatives par analogie au cas associatif.
 
\bigskip

\noindent {\bf Abstract.} We present cogebras of some 
classes of nonassociative algebras whose associator satisfies 
invariance conditions given by the action of the 3-order symmetric group. Amongst these 
algebras we find the well-known Vinberg algebras, the Pre-Lie algebras, 
the Lie-admissible algebras and the 3-power associative algebras. We study properties
between these coalgebras and algebras which generalize the associatif case.

\bigskip

\section{D\'efinitions et exemples}

Dans tout ce travail,
$R$ d\'esigne un anneau commutatif unitaire. Soit $M$ un $R$-module 
et $\Delta$ une comultiplication sur $M$, c'est-\`a-dire une $R$-application lin\'eaire :
$$\Delta: M \rightarrow M\otimes M.$$
Pour tout $R$-module $M,$ nous noterons $\tau:M\otimes M \rightarrow M\otimes M$ l'application lin\'eaire d\'efinie
par $\tau(x \otimes y)=y \otimes x.$ Soit $\Sigma_3$ le groupe sym\'etrique d'ordre 3; notons 
$c_1$ et $c_2$ les $3$-cycles de $\Sigma_3$ et $\tau_{ij}$ la transposition \'echangeant $i$ et $j$.

\begin{defi}
Le couple $(M,\Delta)$ est une cog\`ebre Lie-admissible si l'application lin\'eaire 
$$\Delta_L: M \rightarrow M\otimes M$$
d\'efinie par $\Delta_L=\Delta-\tau \circ \Delta$ est une comultiplication de cog\`ebre de Lie.
\end{defi}
Ainsi $(M,\Delta)$ est une cog\`ebre Lie-admissible si et seulement si $(M,\Delta_L)$ est une cog\`ebre de Lie. 

\noindent Soit $\tilde{A}:End(A,A^{\otimes_2})\rightarrow End(A,A^{\otimes_3})$  le coassociateur donn\'e par 
\begin{eqnarray}
\label{1}
\tilde{A}(\Delta)& = & (\Delta \otimes Id) \circ \Delta -(Id \otimes \Delta) \circ \Delta.
\end{eqnarray}
Pour tout $\sigma \in \Sigma_3$, $\Phi^M_\sigma: M^{\otimes 3} \rightarrow M^{\otimes 3} $ est l'application 
lin\'eaire d\'efinie par 
$$ \Phi^M_\sigma(x_1 \otimes x_2 \otimes x_3)= x_{\sigma^{-1}(1)} \otimes x_{\sigma^{-1}(2)} \otimes x_{\sigma^{-1}(3)}.$$ 
Avec ces notations,
l'homomorphisme $\Delta_L$ v\'erifie
$$
\left\{
\begin{array}{l}
\tau \circ \Delta_L=-\Delta_L \\
\tilde{A}(\Delta_L)+\Phi_{c_1}^M \circ \tilde{A}(\Delta_L)+\Phi_{c_2}^M \circ \tilde{A}(\Delta_L)=0 
\end{array}
\right.
$$

\medskip

\begin{proposition}
Une comultiplication $\Delta$ sur $M$ est une comultiplication Lie-admissible si et seulement si $\Delta$ v\'erifie
\begin{eqnarray}
\label{2}
\sum_{\sigma \in \Sigma_3} (-1)^{\epsilon(\sigma)}\Phi^M_\sigma \circ \tilde{A}(\Delta )=0 
\end{eqnarray}
o\`u $\epsilon(\sigma)$ d\'esigne la signature de la permutation $\sigma$.
\end{proposition}

\noindent {\bf D\'emonstration.} Ceci de d\'eduit directement de l'\'equation~(\ref{1}) sachant que 
$$
\begin{array}{lll}
\tilde{A}(\tau \circ \Delta)& = &((\tau \circ \Delta)\otimes Id)
\circ \tau \circ \Delta  -(Id \otimes (\tau \circ \Delta))\circ (\tau \circ \Delta)\\
& = & -\Phi^M_{ \tau_{13} } \circ \tilde{A}(\Delta).
\end{array}
$$

\bigskip

\noindent {\bf Exemples.}
\begin{itemize}
\item Toute cog\`ebre coassociative est une cog\`ebre Lie-admissible.
\item La comultiplication d'une cog\`ebre pr\'e-Lie $(M,\Delta)$ v\'erifie 
$$ \tilde{A}(\Delta)-\Phi^M_{\tau_{12}}\circ \tilde{A}(\Delta)=0.$$ 
Elle v\'erifie donc aussi l'identit\'e~(\ref{2}) et toute alg\`ebre pr\'e-Lie est Lie-admissible. 
\end{itemize}

\noindent Dans le paragraphe suivant nous g\'en\'eralisons ce type d'exemples.

\section{Les $G_i$-cog\`ebres}

\subsection{Les alg\`ebres $G_i$-associatives}

Notons $G_1=\{Id\},G_2=\{Id,\tau_{12}\},G_3=\{Id,\tau_{23}\},  G_4=\{Id,\tau_{13}\},G_5=\{Id,c_1,c_2\}$ et $G_6=\Sigma_3,$ 
les sous-groupes de $\Sigma_3$. Pour chacun de ces sous-groupes, on d\'efinit  l'application lin\'eaire
$\Phi^M_{G_i}$ de $M^{\otimes 3}$ dans lui-m\^eme par:
$$\Phi^M_{G_i}=\sum_{\sigma \in G_i} (-1)^{\epsilon(\sigma)} \Phi^M_{\sigma}.$$
Soit $\mu:M \otimes M \rightarrow M$ une multiplication de $M$ et notons 
$A:End(M^{\otimes 2},M) \rightarrow End(M^{\otimes 3},M)$        
l'associateur de $\mu$:
$$A(\mu)=\mu \circ (\mu \otimes Id)-\mu \circ (Id \otimes \mu ).$$
Dans \cite{G.R} et \cite{R} nous avons d\'efini les alg\`ebres $G_i$-associatives. Nous dirons que la multiplication $\mu$ est 
$G_i$-associative si son associateur $A(\mu)$ v\'erifie 
\begin{eqnarray}
\label{G_i}
A(\mu) \circ \Phi^M_{G_i}=0.
\end{eqnarray}
En particulier les alg\`ebres $G_1$-associatives sont les alg\`ebres associatives, 
les alg\`ebres $G_2$-associatives sont les alg\`ebres de Vinberg,
les alg\`ebres $G_3$-associatives sont les alg\`ebres pr\'e-Lie,
et les alg\`ebres $G_6$-associatives, c'est-\`a-dire $\Sigma_3$-associatives, sont les alg\`ebres Lie-admissibles.

\medskip

\noindent {\bf Remarque.} Dans \cite{G.R2}, nous avons \'etendu ces relations de non-associativit\'e 
lorsque $R=\K$ est un corps commutatif de caract\'eristique diff\'erente de $2$ et $3$  
en consid\'erant toutes
les relations du type
$$A(\mu) \circ \Phi_v^M=0$$
o\`u $v=\sum_i a_i \sigma_i$  est un vecteur non nul de la $\K$s-alg\`ebre 
$\K[\Sigma_3]$ du groupe $\Sigma_3$ et  
$$\Phi _v^M=\sum_i a_i \Phi^M_{\sigma_i}.$$
Notons $\mathcal{O}(v)$ l'orbite de $v$ dans l'action \`a droite
$$ \Sigma_3 \otimes \K[\Sigma_3] \rightarrow \K[\Sigma_3] $$
et $F_v$ le sous espace de $\K[\Sigma_3]$ engendr\'e par $\mathcal{O}(v).$ C'est un sous espace invariant.

\noindent R\'eciproquement \'etant donn\'e un sous espace invariant $F$ de $\K[\Sigma_3]$, il existe $v \in \K[\Sigma_3]$ 
tel que $F=F_v$. La classification des sous espaces invariants de $\K[\Sigma_3]$ donne une classification des relations 
de non associativit\'e (voir \cite{G.R2}). En particulier, les sous espaces irr\'eductibles de dimension 1 
correspondent aux relations des alg\` ebres Lie-admissibles et des alg\` ebres \`a  puissance $3$-associative \cite{Al}.

\subsection{Les $G_i^!$-alg\`ebres} 

Notons $G_i-\mathcal{A}ss$ l'op\'erade quadratique associ\'ee aux alg\`ebres $G_i$-associatives. 
Dans \cite{G.R} et \cite{M.R}, on montre que ces op\'erades v\'erifient la dualit\'e de Koszul seulement pour 
$i=1,2,3,6$. Soit $G_i^!-\mathcal{A}ss$ l'op\'erade duale. Les alg\`ebres correspondantes sont associatives et v\'erifient
les relations suivantes:

- pour $i=2$ : $x_1.x_2.x_3=x_2.x_1.x_3,$

- pour $i=3$ : $x_1.x_2.x_3=x_1.x_3.x_2,$

- pour $i= 4$ : $x_1.x_2.x_3=x_3.x_2.x_1,$

- pour $i=5$ : $x_1.x_2.x_3=x_2.x_3.x_1=x_3.x_1.x_2,$

- pour $i=6$ : $x_1.x_2.x_3=x_{\sigma(1)}.x_{\sigma(2)}.x_{\sigma(3)} $ pour tous $x_1,x_2,x_3 \in M$ et
$\sigma \in  \Sigma _3$.

\subsection{Les $G_i$-cog\`ebres}

En dualisant la formule~(\ref{G_i}) on obtient la notion de $G_i$-cog\`ebre.

\begin{definition}
Une  $G_i$-cog\`ebre est un $R$-module $M$ muni d'une comultiplication $\Delta$ v\'erifiant
$$\sum_{\sigma \in G_i} (-1)^{\epsilon(\sigma)}\Phi^M_\sigma \circ \tilde{A}(\Delta)=0.$$
\end{definition}

Soit $V$ le vecteur de $R[\Sigma_3]$ d\'efini par $V=Id-\tau_{12}-\tau_{23}-\tau_{13}+c_1+c_2.$ Les cog\`ebres
Lie-admissibles sont d\'efinies par la relation: 
$$\sum_{\sigma \in \Sigma_3} (-1)^{\epsilon(\sigma)}\Phi^M_\sigma \circ \tilde{A}(\Delta)=
\Phi^M_V \circ \tilde{A}(\Delta)=0.$$

Pour tout $v_i=\sum_{\sigma \in G_i} (-1)^{\epsilon(\sigma)} \sigma$ on a $V \in F_{v_i}$.
Ceci implique 

\begin{proposition}
Toute $G_i$-cog\`ebre est une cog\`ebre Lie-admissible.
\end{proposition}

\section{Propri\'et\'es des $G_i$-cog\`ebres}

\subsection{Sur l'espace dual d'une $G_i$-cog\`ebre}

\begin{proposition}
Soit $(M,\Delta)$ une $G_i$-cog\`ebre. Alors le dual $M^{*}$ du $R$-module $M$ est naturellement muni 
d'une structure d'alg\`ebre $G_i$-associative.
\end{proposition}

\noindent {\it D\'emonstration.} 
Pour tout entier naturel $n$ et tout $R$-modules $E$ et $F$, on note 
$$\lambda _n:Hom(E,F)^{\otimes n} \rightarrow Hom(E^{\otimes n},F^{\otimes n})$$
le plongement naturel
$$ \lambda_n(f_1 \otimes ... \otimes f_n)(x_1 \otimes ... \otimes x_n)=f_1(x_1) \otimes ... \otimes f_n(x_n).$$
Ceci \'etant, munissons $M^{*}$ de la multiplication
$$\mu (f_1 \otimes f_2)=\mu_R \circ \lambda_2(f_1 \otimes f_2) \circ \Delta,$$
o\`u $\mu_R$ est la multiplication de $R$. Alors
$$ \mu  \circ (\mu \otimes Id)(f_1 \otimes f_2 \otimes f_3)=
\mu_R \circ (\mu_R \otimes Id ) \circ \lambda_3(f_1 \otimes f_2 \otimes f_3) \circ (\Delta \otimes Id) \circ \Delta,$$
ce qui implique: 
$$A(\mu )(f_1 \otimes f_2 \otimes f_3)=\mu_R \circ (\mu_R \otimes Id) \circ \lambda_3 
(f_1 \otimes f_2 \otimes f_3) \circ (\Delta \otimes Id-Id \otimes \Delta )\circ \Delta.$$
Ainsi 
$\sum_{\sigma \in G_i} (-1)^{\epsilon(\sigma)}A(\mu) \circ \Phi_\sigma^{M^{*}}(f_1 \otimes f_2 \otimes f_3)$

\smallskip

$\begin{array}{l}
\quad  = \sum_{\sigma \in G_i} (-1)^{\epsilon(\sigma)}  \mu_R \circ (\mu_R \otimes Id) \circ 
\lambda_3(\Phi^{M^*}_\sigma (f_1 \otimes f_2 \otimes f_3) ) \circ (\Delta \otimes Id-id \otimes \Delta) \circ \Delta \\
\smallskip \\
 \quad =
\sum_{\sigma \in G_i} (-1)^{\epsilon(\sigma)}\mu_R \circ (\mu_R \otimes Id) \circ \lambda_3(f_1 \otimes f_2 \otimes f_3) 
\circ \Phi^{M}_\sigma \circ (\Delta \otimes Id-id \otimes \Delta) \circ \Delta \\
\smallskip \\
 \quad =
\mu_R \circ (\mu_R \otimes Id) \circ \lambda_3(f_1 \otimes f_2 \otimes f_3)  \circ 
\sum_{\sigma \in G_i} (-1)^{\epsilon(\sigma)} \Phi^{M}_\sigma  \circ \tilde{A}(\Delta)=0. \qquad \ \rule{0.5em}{0.5em}
\end{array}$

\begin{proposition}
Soit $M$ un $R$-module libre de type fini muni d'une structure d'alg\`ebre $G_i$-associative. Alors le dual $M^{*}$
est une $G_i$-cog\`ebre.
\end{proposition}
{\it D\'emonstration.} Soit $M $ une alg\`ebre $G_i$-associative  de dimension finie et soit
$\{{e_i},{i=1,...,n}\}$ une base $M $. Si
$\{f_i\}$ est la base duale, alors $\{f_i \otimes f_j\}$ est une base de
$M^* \otimes M^* $. 
Le coproduit $\Delta $ sur $M^*$ est d\'efini par 
$$\Delta (f)=\sum _{i,j} f(\mu (e_i \otimes e_j)) f_i \otimes f_j.$$
En particulier
$$\Delta (f_k) = \sum _{i,j} C_{ij}^{k}f_i \otimes f_j$$
o\`u $C_{ij}^k$ sont les constantes de structure  de $\mu$ relatives \`a la base $\{{e_i}\}$. 
Ainsi $\Delta$ est la comultiplication d'une  $G_i$-cog\`ebre. \ \rule{0.5em}{0.5em}

\medskip

\subsection{Produit de convolution}

Rappelons que si $(\mathcal{A},\mu)$ est une $R$-alg\`ebre associative et   $(C,\Delta)$ une 
$R$-cog\`ebre coassociative (i.e. une $G_1$-cog\`ebre) alors le produit
de convolution 
$$
f \star g = \mu \circ \lambda_2 \circ  (f\otimes g)\circ \Delta 
$$
munit $Hom(C,\mathcal{A})$ d'une structure d'alg\`ebre associative. Ce r\'esultat s'\'etend au cas des 
$G_i$-alg\`ebres et cog\`ebres apr\`es avoir remarqu\'e que la classe des  cog\`ebres coassociatives est la classe
des  $G_1^{!}$-cog\`ebres.

\subsubsection{$G_i^!$-cog\`ebres}
Dualisons la relation de d\'efinition d'une $G_i^!$-alg\`ebre introduite au paragraphe pr\'ec\'edent.
\begin{defi}
Pour tout $i$, $1 \leq i \leq 6$, une $G_i^!$-cog\`ebre  est une cog\`ebre coassociative, c'est-\`a-dire une $G_1$-cog\`ebre,
v\'erifiant :
$$\Phi_{u_i} \circ (Id\otimes \Delta )\circ \Delta =(Id\otimes \Delta )\circ \Delta $$
o\`u $u_i=\sum _{\sigma \in G_i}\sigma^{-1}.$
\end{defi}

\subsubsection{Structure d'alg\`ebre $G_i$-associative sur $Hom(C,\mathcal{A})$}

\begin{proposition}
Soit $(\mathcal{A},\mu)$ une alg\`ebre  $G_i$-associative et $(C,\Delta )$
une $G_i^!$-cog\`ebre. Alors  $(Hom(C,\mathcal{A}), \star)$ est munie d'une structure d'alg\`ebre 
$G_i$-associative o\`u $\star $ est le produit de convolution:
$$f\star g=\mu \circ \lambda _2(f\otimes g)\circ \Delta .$$
\end{proposition}

\noindent {\it D\'emonstration.}  Calculons l'associateur $A_*$ du produit de convolution.
$$
\begin{array}{ll}
(f_1 \star f_2) \star f_3 & = \mu \ \circ \lambda _2 ((f_1 \star f_2)\otimes f_3 )\,  \circ \,
\Delta
\\ 
& \vspace{0.1cm}\\
& =\mu  \, \circ \, \lambda _2 ((\mu\circ \lambda _2 (f_1 \otimes f_2)\circ  \Delta) \otimes f_3 ) \,  \circ \,
\Delta  \\
& \vspace{0.1cm}\\
& =\mu  \, \circ \,  (\mu \otimes Id)\circ \lambda _3(f_1 \otimes f_2 \otimes
f_3)\circ (\Delta \otimes Id)\, \circ \, \Delta.
\end{array}
$$
On a donc 
$$
\begin{array}{ll}
A_{\star }(f_1\otimes f_2\otimes f_3)  &=\mu  \, \circ \,  (\mu \otimes Id)\circ \lambda _3(f_1
 \otimes f_2 \otimes f_3)\circ (\Delta \otimes Id)\, \circ \, \Delta\\
& \vspace{0.1cm}\\
& - \mu  \, \circ \,  (Id \otimes \mu )\circ \lambda _3(f_1 \otimes f_2 \otimes
f_3)\circ (Id\otimes \Delta )\, \circ \, \Delta.
\end{array}
$$
Alors
$$
\begin{array}{llllllllllllllllllll}
\sum_{\sigma \in G_i}(-1)^{\epsilon (\sigma )}A_{\star } \circ \Phi_{\sigma }^{Hom(C,A)}(f_1\otimes f_2\otimes f_3) 
&&&&&&&&&&&&&&&&&&&
\end{array}
$$
$$\begin{array}{lll}
& =&\mu  \, \circ \,  (\mu \otimes Id)\circ 
(\sum_{\sigma \in G_i}\lambda _3( \Phi_{\sigma }^{Hom(C,A)}(f_1
 \otimes f_2 \otimes f_3)))\circ (\Delta \otimes Id)\, \circ \, \Delta\\
& &\vspace{0.05cm}\\
&& - \mu  \, \circ \,  (Id \otimes \mu )\circ (\sum_{\sigma \in G_i} \lambda _3(
\Phi_{\sigma }^{Hom(C,A)}(f_1 \otimes f_2 \otimes
f_3)))\circ (Id\otimes \Delta )\, \circ \, \Delta.
\end{array}
$$
Mais
$$\lambda _3(
\Phi_{\sigma }^{Hom(C,A)}(f_1 \otimes f_2 \otimes
f_3))=\Phi_{\sigma }^{A}\circ \lambda _3(f_1 \otimes f_2 \otimes
f_3)\circ \Phi_{\sigma^{-1} }^{C}.$$
Ceci implique 
$$
\begin{array}{llllllllllllllllllll}
\sum_{\sigma \in G_i}(-1)^{\epsilon (\sigma )}A_{\star } \circ \Phi_{\sigma }^{Hom(C,A)}(f_1\otimes f_2\otimes f_3) 
&&&&&&&&&&&&&&&&&&&
\end{array}
$$
$$\begin{array}{lll}
& =&\mu  \, \circ \,  (\mu \otimes Id)\circ 
(\sum_{\sigma \in G_i}
\Phi_{\sigma }^{A}\circ \lambda _3(f_1 \otimes f_2 \otimes
f_3))\circ \Phi_{\sigma^{-1} }^{C}\circ 
(\Delta \otimes Id)\, \circ \, \Delta\\
&&\vspace{0.05cm}\\
&& - \mu  \, \circ \,  (Id \otimes \mu )\circ (\sum_{\sigma \in G_i} 
\Phi_{\sigma }^{A}\circ \lambda _3(f_1 \otimes f_2 \otimes
f_3))\circ \Phi_{\sigma^{-1} }^{C}\circ (Id\otimes \Delta )\, \circ \, \Delta.
\end{array}
$$
\smallskip

\noindent Comme $\Delta $ est coassociative 
$$(\Delta \otimes Id)\, \circ \, \Delta=(Id\otimes \Delta )\, \circ \, \Delta$$
et la  structure de $G_i^!$-cog\`ebre  implique
$$\Phi_{u_i }^{C}\circ (Id\otimes \Delta )\, \circ \, \Delta=
\Phi_{u_i }^{C}\circ 
(\Delta \otimes Id)\, \circ \, \Delta=
(\Delta \otimes Id)\, \circ \, \Delta.$$
Alors
$$
\begin{array}{llllllllllllllllllll}
\sum_{\sigma \in G_i}(-1)^{\epsilon (\sigma )}A_{\star } \circ \Phi_{\sigma }^{Hom(C,A)}(f_1\otimes f_2\otimes f_3) 
&&&&&&&&&&&&&&&&&&&
\end{array}
$$
$$
\begin{array}{lll}
& =&\mu  \, \circ \,  (\mu \otimes Id)\circ 
(\sum_{\sigma \in G_i}
\Phi_{\sigma }^{A}\circ \lambda _3(f_1 \otimes f_2 \otimes
f_3))\circ (\Delta \otimes Id)\, \circ \, \Delta\\
&&\vspace{0.05cm}\\
&& - \mu  \, \circ \,  (Id \otimes \mu )\circ (\sum_{\sigma \in G_i} 
\Phi_{\sigma }^{A}\circ \lambda _3(f_1 \otimes f_2 \otimes
f_3))\circ (\Delta \otimes Id)\, \circ \, \Delta\\
& &\smallskip\\
&=&\sum_{\sigma \in G_i}A_{\mu }\circ \Phi_{\sigma }^{A}\circ \lambda _3(f_1 \otimes f_2 \otimes
f_3)\circ (\Delta \otimes Id)\, \circ \, \Delta \\
&&\smallskip \\
&=&0.
\end{array}
$$
\smallskip

\noindent D'o\`u la proposition. \ \rule{0.5em}{0.5em}

\bigskip

\noindent {\bf Remarque:  Big\`ebres Lie-admissibles.} Une big\`ebre Lie-admissible est d\'efinie par la donn\'ee
d'un triplet $(\mathcal{A},\mu,\Delta)$ o\`u  $(\mathcal{A},\mu)$ est une alg\`ebre Lie-admissible, 
$(\mathcal{A},\Delta)$ une cog\`ebre Lie-admissible et d'une relation de compatibilit\'e
liant $\Delta$ et $\mu$ c'est \`a dire on d\'efinit $\Delta \circ \mu$ de sorte que
$$\Delta \circ A(\mu) \circ \Phi_{G_6}^{\mathcal{A}}=0. $$
Nous ne demandons pas \`a l'alg\`ebre d'\^etre unitaire ou \`a la cog\`ebre d'\^etre counitaire.
Parmi ces big\`ebres Lie-admissibles, on aura la classe des $G_i$-big\`ebres. A titre d'exemple, une relation 
de compatibilit\'e pour
les pr\'e-Lie big\`ebres (c'est \`a dire $G_3$-big\`ebres)
est donn\'ee par
$$\Delta \circ \mu=(Id \otimes \mu) \circ (\Delta \otimes Id)+(\mu \otimes Id) \circ 
\Phi_{\tau_{23}}^{\mathcal{A}} \circ(\Delta \otimes Id).$$

\end{document}